\documentclass[12pt]{article}

\def \PTP{Prog.~Theor.~Phys.~}
\def\CRAS{C.~R.~Acad.~Sc.~Paris}
\def \jmax{J}
\def \PII    {{\rm P2}}

\def \LHS{l.h.s.~}
\def \LHS{lhs~}
\def \D {\hbox{d}}
\def \mod#1{\vert #1 \vert}
\def \sech{\mathop{\rm sech}\nolimits}
\def \coth{\mathop{\rm coth}\nolimits}
\def \Degree{\mathop{\rm deg}\nolimits}
\def \Order{\mathop{\rm order}\nolimits}

\begin{document}

\title{Elliptic general analytic solutions\footnote
{
To appear, Studies in Applied Mathematics. Corresponding author RC,
fax +33--169088786.
}
}

\author{Robert Conte\dag\ and M.~Musette\ddag\
{}\\
\\ \dag
LRC MESO,
\\ \'Ecole normale sup\'erieure de Cachan (CMLA) et CEA--DAM
\\ 61, avenue du Pr\'esident Wilson,
\\ F--94235 Cachan Cedex, France.
\\and
\\ Department of Mathematics, The University of Hong Kong
\\ Pokfulam, Hong Kong
\\ E-mail:  Robert.Conte@cea.fr
{}\\
\\ \ddag Dienst Theoretische Natuurkunde,
Vrije Universiteit Brussel,
Pleinlaan 2,
\\
\noindent
B-1050 Bruxelles
\\ E-mail: MMusette@vub.ac.be
}

\maketitle

{\vglue -10.0 truemm}
{\vskip -10.0 truemm}

\begin{abstract}
In order to find analytically 
the travelling waves of 
partially integrable autonomous
nonlinear partial differential equations,
many methods have been proposed over the ages:
``projective Riccati method'',
``tanh-method'', 
``exponential method'', 
``Jacobi expansion method'', 
``new ...'', etc.
The common default to all these ``truncation methods''
is to only provide some solutions, not all of them.
By implementing three classical results of Briot, Bouquet and Poincar\'e,
we present an algorithm able to provide in closed form
\textit{all} those travelling waves which are elliptic
or degenerate elliptic, i.e.~rational in one exponential or rational.
Our examples here include  
the Kuramoto-Siva\-shinsky equation
and
the cubic and quintic complex Ginzburg-Landau equations.

\end{abstract}



\noindent \textit{PACS 2001} 
  02.30.Jr,
  02.30.Ik,
  42.65.Wi.


\baselineskip=12truept

\tableofcontents

\section{Introduction} 
\label{sectionIntro}

In physical language,
given an autonomous 
nonlinear partial differential equation (PDE)
$E(v,v_x,v_t,\dots,v_{n x, m t})=0$,
one calls 
\textbf{travelling wave} any solution of the reduction
   $(v,x,t) \to \xi=x-ct, v(x,t)=u(\xi)$, with $c$ constant,
and \textbf{solitary wave} a travelling wave
which obeys the decay condition that $u(\xi)$ should have constant
limits, possibly different, as $\xi \to \pm \infty$.
The problem addressed here is
to find all the travelling waves of a given PDE, in closed form,
or more generally all the solutions of similarity reductions
of a given PDE to some autonomous ODE.
Two cases naturally arise.
If the PDE is ``integrable'',
for instance in the sense of the inverse spectral transform
\cite{AblowitzClarkson},
there exist powerful methods to solve this problem completely; 
we will discard this case.
If the PDE is partially integrable
(i.e.~fails at least one of the integrability criteria),
there exist several nonperturbative methods
(mostly ``truncations'' \cite{WTC})
able to find \textit{some} travelling waves
but generically unable to find \textit{all} of them.

In a more mathematical language,
let us denote $N$ the differential order of the autonomous
ordinary differential equation (ODE)
\begin{eqnarray}
& &
E(u,u',\dots,u^{(N)})=0,\ '=\frac{\D}{\D \xi},
\label{eqODE}
\end{eqnarray}
which is the travelling wave reduction of the considered PDE.
If the ODE has the Painlev\'e property \cite{CMBook}
(i.e.~no ``bad'' singularities),
the problem is to find its general solution
(depending on $N$ arbitrary constants) in closed form; we discard this case.
If the ODE fails the Painlev\'e test (``partially integrable'' case),
by definition it admits some movable critical singularities
(movable: which depends on the initial conditions;
critical: the solution is multi-valued around it).
Let us define its
\textbf{general analytic solution}
as the $M$-parameter particular solution without movable critical
singularities, with $M$ maximal and of course $M<N$.
The problem is then to find this general analytic solution in closed form.
This is equivalent to find the largest order subequation of the given ODE
with the Painlev\'e property
(subequation: the ODE is a differential consequence of it).
\medskip

The present article, which is self-contained,
gives a synthetic presentation of a method \cite{MC2003,CMBook} 
able to find this general analytic solution
(i.e.~\textbf{all} the single-valued travelling waves)
when it is elliptic or degenerate of elliptic.
This includes the physically relevant solitary waves.

There are many equivalent definitions for elliptic functions.
{}From the point of view of nonlinear differential equations adopted here,
the suitable definition for an \textbf{elliptic function} $f(\xi)$
is any rational function of $\wp(\xi),\wp'(\xi)$,
\begin{eqnarray}
& &
f(\xi)=R(\wp(\xi),\wp'(\xi)),
\label{eqRationalwpwpprime}
\end{eqnarray}
in which $\wp$ is the function of Weierstrass,
defined by the first order ODE
\begin{eqnarray}
& &
{\wp'}^2= 4 \wp^3 - g_2 \wp - g_3
        =4(\wp(\xi)-e_1)(\wp(\xi)-e_2)(\wp(\xi)-e_3),\
\wp'' = 6 \wp^2 - \frac{g_2}{2}.
\end{eqnarray}
Elliptic functions are doubly periodic and have two successive degeneracies,
\begin{itemize}
\item
when one root $e_j$ is double ($g_2^3-27 g_3^2=0$),
degeneracy to simply periodic functions
(i.e.~rational functions of one exponential $e^{k x}$)
according to
\begin{eqnarray}
& &
\forall x,d:\
\wp(x,3 d^2,-d^3)
 = 2 d         - \frac{3 d}{2} \coth^  2  \sqrt{\frac{3 d}{2}} x,
\label{eqwpcoth}
\end{eqnarray}

\item
when the root $e_j$ is triple ($g_2=g_3=0$),
degeneracy to rational functions of $\xi$.

\end{itemize}

In section \ref{section2},
we handle a textbook example of an integrable equation 
in order to make several points.
In section \ref{section3},
we make precise the notion of \textit{general analytic solution}
of a partially integrable ODE.
In section \ref{sectionExisting},
we shortly review the existing methods and point out
their weaknesses.
In section \ref{sectionSubeq},
we present the algorithmic method,
then we process a few selected examples.
Finally, section \ref{sectionExtensions}
examines possible extensions of the method.

\section{Tutorial integrable example} 
\label{section2}

Consider the Korteweg-de Vries (KdV) equation
\begin{eqnarray}
& &
b v_t + v_{xxx} - \frac{6}{a} v v_x=0,\ (a,b) \hbox{ constant},
\label{eqKdV}
\end{eqnarray}
and its travelling wave reduction
\begin{eqnarray}
& &
v=-\frac{a b}{6} c + u,\ \xi=x-c t,\
u''' - \frac{6}{a} u u'=0.
\label{eqKdVreducelliptic}
\end{eqnarray}

\subsection{Mathematical part} 

If one \textit{a priori} knows nothing on (\ref{eqKdVreducelliptic}),
some skill is needed to find the successive first integrals
$g_2$ (arising from the conservative form)
\begin{eqnarray}
& &
u'' - \frac{3}{a} u^2 + a g_2=0,
\label{eqKdVreducellipticwithg2}
\end{eqnarray}
and $g_3$ (integrating factor $u'$)
\begin{eqnarray}
& &
\frac{1}{2} {u'}^2 - \frac{1}{a} u^3 + a g_2 u + 2 a^2 g_3=0.
\label{eqKdVreducellipticOrder1}
\end{eqnarray}
The general solution is then elliptic (doubly periodic)
\begin{eqnarray}
& &
u(\xi)=2 a \wp(\xi-\xi_0,g_2,g_3),
\label{eqKdVreducellipticwp}
\end{eqnarray}
in which the three constants of integration $\xi_0,g_2,g_3$
are allowed to be complex. 

\subsection{Physical part} 

One then prescribes some behaviour at $\xi \to \pm \infty$, e.g.
\begin{eqnarray}
& & {\hskip -8.0 truemm}
\lim_{\xi \to - \infty} u =B_-,\
\lim_{\xi \to + \infty} u =B_+,\ B_\pm \hbox{ real},
\end{eqnarray}
which for $B_-=B_+=B$ makes the doubly periodic function degenerate to
simply periodic (trigonometric)
\begin{eqnarray}
& &
\left\lbrace
\begin{array}{ll}
\displaystyle{
g_2=3 ((B-\kappa)/a)^2,\
g_3=- ((B-\kappa)/a)^3,\
\kappa=-a b c/6.
}\\ \displaystyle{
u=- \frac{a b}{6} c+2 a \left(\frac{4}{3} k^2
  - k^2 \coth^2 k (\xi-\xi_0)\right),\
k^2=\frac{3(B-\kappa)}{2 a},
}
\end{array}
\right.
\end{eqnarray}
in which the constants of integration are $k^2$ (real) and $\xi_0$ (complex).

Finally one requires this solution $u$ to be bounded on the real axis,
by moving all singularities outside $\mathcal{R}$
\begin{eqnarray}
& & {\hskip -11.0 truemm}
\xi_0=\xi_1 + i  \pi/(2 k),\
\xi_1 \in \mathcal{R},\
\coth k (\xi-\xi_0)=\tanh k (\xi-\xi_1),
\end{eqnarray}
to obtain the physical, bell-shaped, solitary wave on the background $B$,
\begin{eqnarray}
& & {\hskip -8.0 truemm}
u=B - 2 a k^2 \sech^2 
     \left(k x-4 k^3 t/b + 6 k (B/a) t/b -k \xi_1\right).
\label{eqsech2KdV}
\end{eqnarray}

\subsection{Lessons from this tutorial example} 

Although it is completely straightforward, 
this example nevertheless teaches us several lessons.
\begin{itemize}    
\item 
Skill was needed to find the first integrals, we want to avoid that.

\item
If possible, it is better to first find the mathematical solution,
because all the physical solutions follow.
In the case of the nonlinear Schr\"odinger equation for instance,
which admits two
(one bright, one dark)
solitary waves depending on the focusing or defocusing situation,
the bright wave is found as above,
and the dark one arises from a factorisation of the first order
subequation analogous to (\ref{eqKdVreducellipticOrder1}).
  
\item
Elliptic functions seem privileged.
We will come back to that later.

\end{itemize}

\section{Partially integrable case}
\label{section3}

By definition,
partially integrable PDEs fail the Painlev\'e test.
Our specific examples will be the following reductions: 
\begin{itemize}    
\item the similarity reduction of the
one-dimensional cubic complex Ginzburg-Landau equation (CGL3)
\begin{eqnarray}
& & {\hskip -10.0 truemm}
i A_t + p A_{xx} + q \mod{A}^2 A - i \gamma A =0,\
p q \gamma \not=0,\
(A,p,q) \in {\mathcal C},\
\gamma  \in {\mathcal R}
\label{eqCGL3}
\\
& &
A(x,t)=\sqrt{M(\xi)} e^{i(\displaystyle{-\omega t + \varphi(\xi)})},\
\xi=x-ct,
\label{eqCGL3red}
\end{eqnarray}
(a generic equation for slowly varying amplitudes,
see the review \cite{vS2003}); 

\item the travelling wave reduction of
the Kuramoto-Sivashinsky equation (KS)
\begin{eqnarray}
& &
v_t + \nu v_{xxxx} + b v_{xxx} + \mu v_{xx} + v v_x =0,\ \nu \not=0,
\label{eqKS}
\\
& &
v(x,t)=c+u(\xi),\ \xi=x-ct,\
 \nu u''' + b u'' + \mu u' + \frac{u^2}{2} + A = 0,\ \nu \not=0,
\label{eqKSODE}
\end{eqnarray}
in which $A$ is an integration constant
(flame on a vertical wall, phase equation of CGL3,
see \cite{MannevilleBook} for some physical background);

\item the similarity reduction of 
the one-dimensional quintic complex Ginzburg-Landau equation (CGL5) \cite{vS2003},
\begin{eqnarray}
& &
i A_t +p A_{xx} +q \mod{A}^2 A +r \mod{A}^4 A -i \gamma A =0,\
p r \not=0,\
\Im(p/r)\not=0,\
\label{eqCGL5}
\\
& &
(A,p,q,r) \in {\mathcal C},\
\gamma  \in {\mathcal R}.
\nonumber
\\
& &
A(x,t)=\sqrt{M(\xi)} e^{i(\displaystyle{-\omega t + \varphi(\xi)})},\
\xi=x-ct.
\label{eqCGL5red}
\end{eqnarray}

\end{itemize}

\subsection{General analytic solution of a partially integrable ODE}

Given an $N$-th order ODE (\ref{eqODE}) which is partially integrable
and therefore displays movable multi-valuedness,
one first needs to count the number of integration constants
which correspond to this multi-valuedness.
Let us take an example.

\subsection{Local representation of the general analytic solution}

The KS PDE (\ref{eqKS}) admits the
travelling wave reduction
(\ref{eqKSODE}).
It has a chaotic behavior \cite{MannevilleBook},
and it depends on two dimensionless parameters,
$b^2/(\mu \nu)$ and $\nu A / \mu^3$.
The assumption of a singular algebraic behaviour
\begin{eqnarray}
& &
\xi \to \xi_0:\ u \sim u_0 \chi^p,\ \chi=\xi - \xi_0,\ u_0 \not=0,\ 
p \notin \mathcal{N},
\end{eqnarray}
generates for $(p,u_0)$ the system
\begin{eqnarray}
& &
p-3 = 2p,\
\nu p(p-1)(p-2) u_0 + (1/2) u_0^2=0,\
u_0 \not=0,
\label{eqKSLeadingOrderEq}
\end{eqnarray}
which admits the unique solution $p=-3,u_0=120 \nu$,
with the common value of the two powers: $q=p-3=2p=-6$.
Near this triple pole, the linearized equation has the indicial equation
\begin{eqnarray}
& &
\lim_{\chi \to 0} \chi^{-j-q} (\nu \partial_\xi^3 + u_0 \chi^p) \chi^{j+p}
\\
& &
= \nu (j-3)(j-4)(j-5) + 120 \nu = \nu (j+1) (j^2 -13 j + 60)
\\
& &
=\nu (j+1) \left(j-\frac{13 + i \sqrt{71}}{2}\right)
           \left(j-\frac{13 - i \sqrt{71}}{2}\right)=0,
\label{eqKSIndicial}
\end{eqnarray}
i.e.~two of the three Fuchs indices are irrational complex numbers.
Computing the next terms after the dominant one yields the Laurent series
\begin{eqnarray}
{\hskip -12.0 truemm}
& &
u^{(0)} = 120 \nu \chi^{-3} - 15b \chi^{-2}
        + \frac{15 (16 \mu \nu - b^2)}{4 \times 19 \nu} \chi^{-1}
        + \frac{13 (4  \mu \nu - b^2) b}{32 \times 19 \nu^2}
        + O(\chi),
\label{eqKSODELaurent}
\end{eqnarray}
which only depends on one integration constant.
Laurent series such as (\ref{eqKSODELaurent}) are convergent, 
as proven by Chazy \cite{ChazyThese}.
By perturbing this solution \cite[pages 338, 384]{ChazyThese} 
\cite{CFP1993},
one obtains the three-parameter solution
\begin{eqnarray}
u(\xi-\xi_0,\varepsilon c_+,\varepsilon c_-)
& = &
 u^{(0)}
\nonumber
\\
& &
+ \varepsilon \chi^{-3}[ 
              c_{-1} \chi^{-1} \hbox{Regular}(\chi)
\nonumber
\\
& &
        + \ \ \ \ \ \ \ c_{+} \chi^{(13+i\sqrt{71})/2} \hbox{Regular}(\chi)
\nonumber
\\
& &
        + \ \ \ \ \ \ \ c_{-} \chi^{(13-i\sqrt{71})/2} \hbox{Regular}(\chi) ]
            + \mathcal{O}(\varepsilon^2)\},
\label{eqKSODELaurentPerturbed}
\end{eqnarray}
which is a local representation valid for $\xi \to \xi_0$
and $\varepsilon \to 0$,
depending on the three arbitrary constants
$\xi_0,\varepsilon c_{+},\varepsilon c_{-}$
(indeed, as shown by Poincar\'e\footnote{ 
With his ``\'equation aux variations'', see Chazy
\cite[page 338 footnote 1]{ChazyThese}.
}
and Darboux \cite{Darboux1883Eaux}, 
the coefficient of $\varepsilon c_{-1}$ in 
(\ref{eqKSODELaurentPerturbed}) 
is the derivative of $u^{(0)}$ with respect 
to its arbitrary constant $\xi_0$,
therefore $c_{-1}$ only represents a perturbation of $\xi_0$
and it can be set to $0$).
The only way for 
(\ref{eqKSODELaurentPerturbed}) to become single-valued
is to cancel the contribution of the two irrational Fuchs indices,
thus restricting (\ref{eqKSODELaurentPerturbed})
to (\ref{eqKSODELaurent}).

The general analytic solution as defined in section \ref{sectionIntro}
depends on $3-2=1$ integration constant,
it is locally represented by the Laurent series (\ref{eqKSODELaurent})
and the problem addressed is to obtain a 
closed form expression for $ u^{(0)}$.

\section{Existing methods (sufficient)}
\label{sectionExisting}
   
Given a nonlinear ODE such as (\ref{eqODE}),
all the methods able to yield single-valued solutions
in closed form
necessarily make use of the structure of movable singularities.
Many ``new'', ``extended'', etc,
methods claim not to do so and to be original,
but this is not true:
no method can avoid the analytic structure of singularities.

Moreover,
in order to be nonperturbative,
i.e.~to only require a finite number of steps,
the integration methods have better use an
elementary ``unit of integration''
which has singularities (such as $\tanh$ or $\sech$)
rather than an entire function (such as $\exp$). 

When nothing is known about the ODE (\ref{eqODE}),
most methods follow the pioneering work of 
Weiss, Tabor and Carnevale \cite{WTC}.
The restriction of their method to ODEs was certainly known 
to the classical authors
involved with explicit integration (Painlev\'e, Gambier, Chazy, Bureau),
although we could not find an explicit, convincing reference.

All these methods proceed as follows:
\begin{enumerate}
\item
Assume a given class of expressions for the general analytic solution,
\item
Check whether there are indeed solutions in that class. 
\end{enumerate}

Various names are used to describe these methods:
truncation method, $\tanh$ method, 
Jacobi expansion method, \dots.

Typical classes of assumed expressions are:
\begin{itemize}
\item
polynomials in $\tanh$
(Weiss \textit{et al.}~\cite{WTC} and followers \cite{Conte1989}),

\item
polynomials in $\wp,\wp'$ 
\cite{FournierSpiegelThual,KudryashovElliptic},

\item
polynomials in $\tanh$ and $\sech$ \cite{JeffreyXu,CM1992,Pickering1993}.

\end{itemize}

We will call all these methods \textbf{sufficient}
because, by construction,
they cannot find any solution outside the given class of expressions.

\subsection{Examples of sufficient methods}
\label{sectionKSTruncation}

The KS ODE (\ref{eqKSODE})
has one family of movable triple poles, 
and the Weierstrass function $\wp$ one family of double poles.
Therefore it is consistent to assume
\begin{eqnarray}
& &
u= c_0 \wp' + c_1 \wp +c_2,\ c_0 \not=0,
\nonumber
\\
& &
E(u) \equiv \nu u''' + b u'' + \mu u' + \frac{u^2}{2} + A = 0,\
\nu \not=0,
\label{eqKSODEbis}
\\
& &
\wp'' = 6 \wp^2 - \frac{g_2}{2},\
{\wp'}^2= 4 \wp^3 - g_2 \wp - g_3.
\nonumber
\end{eqnarray}
The elimination of $\wp''$ and ${\wp'}^2$
transforms the \LHS\ of (\ref{eqKSODEbis}) into
\begin{eqnarray}
& &
E(u) =\sum_{j=0}^{k} \sum_{k=0}^{1} E_{j,k} u^j {u'}^k,
\nonumber
\end{eqnarray}
and one has to solve the six determining equations
\begin{eqnarray}
& &
\left\lbrace
\begin{array}{ll}
\displaystyle{
E_{3,0} \equiv c_0 (120 \nu + 2 c_0) =0,
}
\\
\displaystyle{
E_{1,1} \equiv 12 b c_0 + c_0 c_1 + 12 \nu c_1=0,
}
\\
\displaystyle{
E_{2,0} \equiv 6 \mu c_0 + 6 b c_1 +\frac{c_1^2}{2} =0,
}
\\
\displaystyle{
E_{0,1} \equiv c_0 c_2 + \mu c_1=0,
}
\\
\displaystyle{
E_{1,0} \equiv c_1 c_2 -\frac{1}{2} g_2 c_0^2 - 18 \nu g_2 c_0=0,
}
\\
\displaystyle{
E_{0,0} \equiv A - 12 \nu g_3 c_0
+\frac{1}{2}\left(c_2^2 - b g_2 c_1 - \mu g_2 c_0 -g_3 c_0^2 \right)=0.
}
\end{array}
\right.
\nonumber
\end{eqnarray}
As always in these truncations \cite{WTC},
these equations \textit{must} be solved by decreasing value of 
the singularity degree $3 j + 4 k$, 
for the successive unknowns $c_0,c_1,c_2,g_2,g_3$.
The result is
\cite{FournierSpiegelThual,KudryashovElliptic}
\begin{eqnarray}
& &
\left\lbrace
\begin{array}{ll}
u=\displaystyle{- 60 \nu \wp' - 15 b \wp - \frac{b \mu}{4 \nu}},\
\\
g_2=\displaystyle{\frac{\mu^2}{12 \nu^2}},\
g_3=\displaystyle{\frac{13 \mu^3 + \nu A}{1080 \nu^3}},
b^2=16 \mu \nu,
\end{array}
\right.
\label{eqKSElliptic}
\end{eqnarray}
and this would be the general analytic solution
if $b^2-16 \mu \nu$ were unconstrained.
Removing this constraint is still an open problem.

Instead of considering the Weierstrass equation,
one may consider the Riccati equation with constant coefficients
\begin{eqnarray}
& &
\tau' + \tau^2 + \frac{S}{2}=0,\
S=-\frac{k^2}{2}= \hbox{ constant} \in \mathcal{C},
\label{eqRiccatitau}
\end{eqnarray}
whose singularities are one family of movable simple poles,
and make the similar consistent assumption 
\cite{Conte1988,CM1989}
\begin{eqnarray}
& &
u=\sum_{j=0}^{-p} c_j \tau^{-j-p},\ c_0 \not=0,\ p=-3,
\label{eqKSTruncation_tau}
\end{eqnarray}
with $c_j$ constants to be determined.
The elimination of $\tau'$ from (\ref{eqRiccatitau})
transforms the \LHS\ of (\ref{eqKSODEbis}) into a similar polynomial
\begin{eqnarray}
& &
E(u) = \sum_{j=0}^{-q} E_j \tau^{-j-q}=0,\ q=-6.
\end{eqnarray}
This set of determining equations $E_j=0$
admits six solutions, listed in Table \ref{TableKS},
all represented by
\begin{eqnarray}
u & = &
 120 \nu \tau^3
 - 15 b \tau^2
 +\left(\frac{60}{19} \mu - 30 \nu k^2 - \frac{15 b^2}{4 \times 19 \nu}\right)
 \tau
\nonumber
\\
& &
 + \frac{5}{2} b k^2 - \frac{13 b^3}{32 \times 19 \nu^2}
  + \frac{7 \mu b}{4 \times 19 \nu},\
\tau  =  \frac{k}{2} \tanh \frac{k}{2} (\xi-\xi_0).
\label{eqKSTrigo}
\end{eqnarray}
The solitary waves $b=0$ were found
by Kuramoto and Tsuzuki \cite{KuramotoTsuzuki},
and the three other values of $b^2/(\mu \nu)$ were added  by Kudryashov
\cite{KudryashovKSFourb}.

Currently, no other solution is known to (\ref{eqKSODEbis}).

\begin{table}[h] 
\caption[Kuramoto-Sivashinsky equation. The six trigonometric solutions.]{
The six trigonometric solutions of KS, Eq.~(\ref{eqKSODE}),
with the notation $k^2=-2 S$.
They all have the form (\ref{eqKSTrigo}).
The last line is a degeneracy of the elliptic solution
(\ref{eqKSElliptic}).
}
\vspace{0.2truecm}
\begin{center}
\begin{tabular}{| c | c | c |}
\hline 
$b^2/(\mu\nu)$ & $\nu A/\mu^3$ & $\nu k^2/\mu$
\\ \hline \hline 
$0$ & $-4950/19^3,\ 450/19^3$ & $11/19,\ -1/19$
\\ \hline 
$144/47$ & $-1800/47^3$ & $1/47$
\\ \hline 
$256/73$ & $-4050/73^3$ & $1/73$
\\ \hline 
$16$ & $-18,\ -8$ & $1,\ -1$
\\ \hline 
\end{tabular}
\end{center}
\label{TableKS}
\end{table}

It is easy to build an example making truncation methods fail.
Knowing that rational functions are usually not assumed classes,
such an example is for instance
\begin{eqnarray}
& &
u=\frac{\tanh(\xi-\xi_0)}{2+\tanh^2(\xi-\xi_0)}=
\hbox{outside usual classes},
\\
& &
2 {u'}^2+(24 u^2 - 3) u' + 72 u^4 - 17 u^2 + 1=0.
\label{eq_u_rational}
\end{eqnarray}

\section{Subequation method (necessary)}
\label{sectionSubeq}

\subsection{The class of elliptic functions}

All solutions found in previous examples are elliptic or degenerate.
This quite frequent result for solutions of autonomous nonlinear ODEs
has a simple explanation, dating back to Lazarus Fuchs.
Indeed, this is a classical result (L.~Fuchs, Poincar\'e, Painlev\'e)
that the only autonomous first order algebraic ODEs having a
single-valued general solution are 
algebraic transforms of either the Weierstrass elliptic equation,
see (\ref{eqRationalwpwpprime}),
or the Riccati equation (\ref{eqRiccatitau}).
These first order ODEs are therefore the elementary building blocks
with which to build solutions to autonomous ODEs of higher order.

Conversely,
given an $N$-th order autonomous algebraic ODE,
is it possible to find \textbf{all} solutions in this class
of elliptic or degenerate elliptic functions?
The present section brings a positive answer.

\subsection{Two results of Briot and Bouquet}
\label{sectionResultsBB} 

Given an elliptic function,
its \textbf{elliptic order} is defined as
the number of poles in a period parallelogram,
counting multiplicity of course.
It is equal to the number of zeros.

We will need two classical theorems.

\textit{Theorem} \cite[theorem XVII p.~277]{BriotBouquet}.
Given two elliptic functions $u,v$ with the same periods
of respective elliptic orders $m,n$,
they are linked by an algebraic equation
\begin{eqnarray}
& &
F(u,v) \equiv
 \sum_{k=0}^{m} \sum_{j=0}^{n} a_{j,k} u^j v^k=0,\ 
\label{eqTwoEllFunctions}
\end{eqnarray}
with $\Degree(F,u)=\Order(v)$, $\Degree(F,v)=\Order(u)$.
If in particular $v$ is the derivative of $u$,
the first order ODE obeyed by $u$ takes the precise form
\begin{eqnarray}
& &
F(u,u') \equiv
 \sum_{k=0}^{m} \sum_{j=0}^{2m-2k} a_{j,k} u^j {u'}^k=0,\ a_{0,m}\not=0.
\label{eqsubeqODEOrderOnePP}
\end{eqnarray}

\textit{Theorem} (Briot and Bouquet, Poincar\'e, Painlev\'e).
If a first order $m$-th degree autonomous ODE has the Painlev\'e property,
\begin{itemize}
\item
it must have the form (\ref{eqsubeqODEOrderOnePP}),

\item
its general solution is either elliptic (two periods)
or rational in one exponential $e^{k x}$ (one period)
or rational in $x$ (no period)
(successive degeneracies $g_2^3-27 g_3^2=0$, then $g_2=0$ in
${\wp'}^2= 4 \wp^3 - g_2 \wp - g_3$).

\end{itemize}

\textit{Remark}.
Equation (\ref{eqsubeqODEOrderOnePP}) is invariant under an arbitrary
homographic transformation having constant coefficients,
this is another useful feature of elliptic equations.

\subsection{General method to find all elliptic solutions}
\label{sectionFirstOrderSubequation}

Consider an $N$-th order autonomous algebraic ODE (\ref{eqODE})
admitting at least one Laurent series
\begin{eqnarray}
& &
u=\chi^p \sum_{j=0}^{+\infty} u_j \chi^j,\ \chi=\xi-\xi_0.
\label{eqLaurent}
\end{eqnarray}

If its general analytic solution, as defined section \ref{sectionIntro},
is elliptic or degenerate elliptic,
there exists an algorithm to find it in closed form, which we now present.
It only requires four ingredients: 
the two theorems of Briot and Bouquet stated section \ref{sectionResultsBB},
the existence of at least one Laurent series,
an algorithm of Poincar\'e to be presented soon.
Its input and output are as follows.

\textbf{Input}: an $N$-th order ($N\ge 2$) any degree autonomous algebraic
ODE admitting a Laurent series.

\textbf{Output}: all its elliptic or degenerate elliptic solutions
in closed form.

The successive steps are \cite{MC2003,CMBook}:
\begin{enumerate}
\item
Find the analytic structure of singularities
(e.g., 4 families of simple poles, 2 of double poles).
Deduce the elliptic orders $m,n$ of $u,u'$.

\item
Compute slightly more than $(m+1)^2$ terms in the Laurent series.

\item
Define the first order $m$-th degree
subequation $F(u,u')=0$ 
(it contains at most $(m+1)^2$ coefficients $a_{j,k}$),
\begin{eqnarray}
& &
F(u,u') \equiv
 \sum_{k=0}^{m} \sum_{j=0}^{2m-2k} a_{j,k} u^j {u'}^k=0,\ a_{0,m}\not=0.
\end{eqnarray}

\item
Require each Laurent series (\ref{eqLaurent}) to obey $F(u,u')=0$,
\begin{eqnarray}
& & {\hskip -10.0 truemm}
F \equiv \chi^{m(p-1)} \left(\sum_{j=0}^{\jmax} F_j \chi^j
 + {\mathcal O}(\chi^{\jmax+1})
\right),\
\forall j\ : \ F_j=0.
\label{eqLinearSystemFj}
\end{eqnarray}
and solve this \textbf{linear overdetermined} system for $a_{j,k}$.

\item
Integrate each resulting ODE $F(u,u')=0$.
\end{enumerate}

Several remarks are in order.

\begin{enumerate}
\item
The fourth step generates a \textit{linear},
infinitely overdetermined,
system of equations $F_j=0$ for the unknown coefficients $a_{j,k}$.
It is quite an easy task to solve such a system,
and this is the key advantage of the present algorithm.

\item
In the fifth step,
two cases arise, depending on the genus of the algebraic curve $F=0$.
If the genus is one (nondegenerate elliptic),
there exists an algorithm of Poincar\'e
which builds explicitly the rational function 
$R$, Eq.~(\ref{eqRationalwpwpprime}).
This algorithm has been implemented by Mark van Hoeij \cite{MapleAlgcurves} 
as the package ``algcurves'' of the computer algebra language Maple.
If the genus is zero (degenerate elliptic),
finding $u$ as a rational function of one exponential $e^{k (\xi-\xi_0)}$
or as a rational function of $\xi-\xi_0$ is classical,
and also implemented in ``algcurves''.

\item
In the third and fourth steps, the requirement can be weakened,
at the price of finding less that the general analytic solution.
Consider for instance the artificial ODE
(which actually occurs in the travelling wave of the modified KdV
equation),
\begin{eqnarray}
& &
E(u) \equiv \hbox{some differential consequence of }
a^2 {u'}^2 -(u^2 + b)^2 + c=0,\
\end{eqnarray}
which admits two Laurent series
\begin{eqnarray}
& &
u= \pm a \chi^{-1} + \dots.
\end{eqnarray}
With the assumption 
$m=2$ in step 3 and ``Require each Laurent series'' in step 4,
one will find the Jacobi solution.
But, with the weaker assumption
$m=1$ in step 3 and ``Require one [not two] Laurent series'' in step 4,
one will find the constraint $c=0$ in $E(u)$.
The interested reader can practice on
the ODE admitting, for instance, the rational solution
\begin{eqnarray}
& &
u= c_1(x-x_1)^{-1} +c_2(x-x_2)^{-2} +(x-x_3)^{-3},\
c_j \hbox{ fixed},\ x_j \hbox{ movable},
\end{eqnarray}
for which the assumptions can be $m=6,5,4,3$.

\end{enumerate}

\section{Examples}

\subsection{Tutorial integrable example: KdV}

The ODE (\ref{eqKdVreducelliptic})
presents one movable double pole $u \sim 2 a \chi^{-2}$,
hence the respective elliptic orders $2$ and $3$ for $u$ and $u'$.
Nine terms (in fact five because of parity) 
will prove sufficient in the computation of the unique
Laurent series,
\begin{eqnarray}
& &
u=2 a \chi^{-2} + U_4 \chi^2 + U_6 \chi^4 + \frac{U_4^2}{6 a} \chi^6 + \dots,
\end{eqnarray}
in which $U_4$ and $U_6$ are arbitrary constants.
Assuming in step 3
\begin{eqnarray}
& & {\hskip -10.0 truemm}
F \equiv
              {u'}^2
+ a_{0,1}     {u'}
+ a_{1,1} u   {u'}
+ a_{0,0}
+ a_{1,0} u
+ a_{2,0} u^2
+ a_{3,0} u^3, a_{0,2}=1,
\label{eqsubeqKdV}
\end{eqnarray}
one generates in step 4
the linear overdetermined system (\ref{eqLinearSystemFj}),
\begin{eqnarray}
& &
\left\lbrace
\begin{array}{ll}
\displaystyle{
F_0 \equiv 16 a^2 a_{0,2} + 8 a^3 a_{3,0}=0,
}
\\
\displaystyle{
F_1 \equiv -8 a^2 a_{1,1}=0,
}
\\
\displaystyle{
F_2 \equiv 4 a^2 a_{2,0}=0,
}
\\
\displaystyle{
F_3 \equiv -4 a a_{0,1}=0,
}
\\
\displaystyle{
F_4 \equiv 2 a a_{1,0} -16 a a_{0,2} U_4 +12 a^2 a_{3,0} U_4=0,
}
\\
\displaystyle{
F_5 \equiv 0,
}
\\
\displaystyle{
F_6 \equiv a_{0,0} + 4 a a_{2,0} U_4 -32 a a_{0,2} U_6 +12 a^2 a_{3,0} U_6=0,
}
\\
\displaystyle{
\dots}
\end{array}
\right.
\end{eqnarray}
whose unique solution is
\begin{eqnarray}
& &
{u'}^2 - (2/a) u^3 + 20 U_4 u + 56 a U_6=0.
\label{eqsubeqKdV0}
\end{eqnarray}
The two arbitrary constants $U_4,U_6$ correspond to the two
first integrals of (\ref{eqKdVreducelliptic}),
but no skill is required to find them,
the process is systematic.

In step 5,
the Maple commands would look like

{\tt with(algcurves);}     \hfill                     load the package 

{\tt genus(eq46,u,uprime);}\hfill     will answer ``1''

{\tt Weierstrassform(eq46,u,uprime,wp,wpprime);}
\hfill\break\noindent
the last command yielding the four formulae of the  
birational transformation between $(u,u')$ and $(\wp,\wp')$,
one of them being precisely (\ref{eqKdVreducellipticwp}).

\subsection{Partially integrable example, KS}
\label{sectionKSsubeqFirstOrder}

The Laurent series of (\ref{eqKSODE}) is (\ref{eqKSODELaurent}).

In the first step, 
the unique family of movable triple poles yields
elliptic orders $3$ and $4$ for respectively $u$ and $u'$.
With the normalization $a_{0,3}=1$, the subequation $F(u,u')=0$
contains ten coefficients.
In step 4,
these ten coefficients are first determined 
by solving the Cramer system of ten equations
$F_j=0,j=0:6,8,9,12$.
The remaining infinitely overdetermined nonlinear system for $(\nu,b,\mu,A)$
contains as greatest common divisor (gcd) $b^2-16 \mu \nu$
which defines a first solution
\begin{eqnarray}
& & {\hskip -10.0 truemm}
\frac{b^2}{\mu \nu}=16,\
u_s=u+\frac{3 b^3}{32 \nu^2},
\nonumber
\\
& & {\hskip -10.0 truemm}
\left(u' + \frac{b}{2 \nu} u_s\right)^2
\left(u' - \frac{b}{4 \nu} u_s\right)
+\frac{9}{40 \nu}
\left(u_s^2 + \frac{15 b^6}{1024 \nu^4} + \frac{10 A}{3}\right)^2=0.
\label{eqKSsubeqgenus1}
\end{eqnarray}
After division by this factor,
the remaining nonlinear system for 
$(\nu,b,\mu,A)$ with $b^2-16 \mu \nu \not=0$
admits exactly four solutions
(stopping the series at $j=16$ is enough to obtain the result),
identical to those listed in Table \ref{TableKS}
Section \ref{sectionKSTruncation},
each solution defining a first order third degree subequation,
\begin{eqnarray}
& &
{\hskip -14.0 truemm}
b=0,\
\nonumber
\\
& &
{\hskip -14.0 truemm}
\left(u' + \frac{180 \mu^2}{19^2 \nu}\right)^2
\left(u' - \frac{360 \mu^2}{19^2 \nu}\right)
+\frac{9}{40 \nu}
\left(u^2 + \frac{30 \mu}{19} u' - \frac{30^2 \mu^3}{19^2 \nu}\right)^2=0,\
\label{eqKSsubeqgenus0first}
\\
& &
{\hskip -14.0 truemm}
b=0,\
{u'}^3
+\frac{9}{40 \nu}
\left(u^2 + \frac{30 \mu}{19} u' + \frac{30^2 \mu^3}{19^3 \nu}\right)^2=0,\
\\
& &
{\hskip -14.0 truemm}
\frac{b^2}{\mu \nu}=\frac{144}{47},\
u_s=u-\frac{5 b^3}{144 \nu^2},\
\left(u' + \frac{b}{4 \nu} u_s\right)^3+\frac{9}{40 \nu} u_s^4=0,\
\\
& &
{\hskip -14.0 truemm}
\frac{b^2}{\mu \nu}=\frac{256}{73},\
u_s=u-\frac{45 b^3}{2048 \nu^2},\
\nonumber
\\
& &
{\hskip -14.0 truemm}
\left(u' + \frac{b}{8 \nu} u_s\right)^2
\left(u' + \frac{b}{2 \nu} u_s\right)
+\frac{9}{40 \nu}
\left(u_s^2+\frac{5 b^3}{1024 \nu^2}u_s + \frac{5 b^2}{128 \nu}u'\right)^2=0.
\label{eqKSsubeqgenus0last}
\end{eqnarray}

In order to integrate the two sets of subequations
(\ref{eqKSsubeqgenus1}),
(\ref{eqKSsubeqgenus0first})--(\ref{eqKSsubeqgenus0last}),
one must first compute their genus,
which is one for (\ref{eqKSsubeqgenus1}),
    and zero for (\ref{eqKSsubeqgenus0first})--(\ref{eqKSsubeqgenus0last}).
Therefore (\ref{eqKSsubeqgenus1}) has the elliptic general solution
(\ref{eqKSElliptic}).

As to the general solution of the four others
(\ref{eqKSsubeqgenus0first})--(\ref{eqKSsubeqgenus0last}),
this is the third degree polynomial (\ref{eqKSTrigo})
in $\tanh k (\xi-\xi_0)/2$
already found by the one-family truncation method.

\textit{Remark}.
Canceling the gcd $b^2-16 \mu \nu$ is equivalent to cancel the residue
in the Laurent series (\ref{eqKSODELaurent}),
therefore KS admits no elliptic solution other than 
(\ref{eqKSElliptic}) \cite{Hone2005}.
Indeed, this is a property of elliptic functions that, 
inside a period parallelogram,
the sum of the residues is necessarily zero.
Conversely,
by implementing the condition that the sum of the residues of 
any differential polynomial of $u$ vanishes,
Hone \cite{Hone2005}
has found an elegant method to isolate all the nondegenerate elliptic 
(genus one) solutions,
proving in particular that CGL3 has no travelling wave which is
nondegenerate elliptic.
This method cannot however also yield the degenerate elliptic solutions.

\subsection{Partially integrable example, CGL5}
\label{sectionSubeqCGL5}

As example of a new solution found by the method of section
\ref{sectionFirstOrderSubequation},
one can quote the unique elliptic (nondegenerate, i.e.~genus one) 
travelling wave (\ref{eqCGL5red}) of the CGL5 equation (\ref{eqCGL5})
\cite{VernovCGL5},
\begin{eqnarray}
& &
\left\lbrace
\begin{array}{ll}
\displaystyle{
M=c_0 \frac{4 \wp^2(\xi) - g_2}{4 \wp^2(\xi) + g_2},\
\varphi'=\frac{c s_r}{2} 
+ \left(-\frac{3 g_2}{4 \wp(2 \xi)} \right)^{1/2},
}\\ \displaystyle{
\hbox{notation }
s_r - i s_i \equiv p^{-1},\
g_r \equiv s_r \gamma + s_i \omega,\
e_i \equiv \Im(r/p),\
}\\ \displaystyle{
\hbox{constraints }
q=0,\
\Re(r/p)=0,\
c s_i=0,\
s_r \omega - s_i \gamma + (c s_r/2)^2=0,\
}\\ \displaystyle{
\hbox{values of } c_0,g_2,g_3:\
c_0^2=\frac{4 g_r}{3 e_i},\
g_2=-\frac{g_r^2}{27},\
g_3=0.
}
\end{array}
\right.
\label{eqCGL5Vernovsol}
\end{eqnarray}
This solution,
obtained by combining the present method with the conditions 
on residues \cite{Hone2005},
exists at the price of five real constraints
among the coefficients $p,q,r,\gamma,\omega,c$ of the differential system
for $(M,\varphi')$.

In this example, care should be taken that the two elliptic subequations,
\begin{eqnarray}
& &
\left\lbrace
\begin{array}{ll}
\displaystyle{
3^4 e_i {M'}^4 - M^2 (3 e_i M^2 - 4 g_r)^3=0,
}\\ \displaystyle{
9 {\psi'}^2 -12 \psi^4 - g_r^2=0,\ \psi=\varphi' - \frac{c s_r}{2},
}
\end{array}
\right.
\end{eqnarray}
respectively of a canonical Briot-Bouquet type 
and of a Jacobi type,
involve elliptic functions with \textit{different} first arguments, 
see $\wp(2 \xi)$ term in (\ref{eqCGL5Vernovsol}).

\subsection{Results for various partially integrable PDEs}
\label{sectionResults}

For the various PDEs considered in the text,
and recalled below for convenience,
\begin{eqnarray}
& &
\hbox{KdV }:\ 
u''' - (6/a) u u'=0,
\nonumber
\\
& &
\hbox{KS }:\ 
 \nu u''' + b u'' + \mu u' + u^2/2 + A = 0,\ \nu \not=0,
\nonumber 
\\
& &
\hbox{CGL3 }:\ 
i A_t + p A_{xx} + q \mod{A}^2 A - i \gamma A =0,\
p q \not=0,\
\Im(q/p) \not=0,
\nonumber
\\
& &
\hbox{CGL5 }:\ 
i A_t +p A_{xx} +q \mod{A}^2 A +r \mod{A}^4 A -i \gamma A =0,\
p r \not=0,\
\Im(p/r)\not=0,\
\nonumber
\\
& &
\hbox{2 CGL3 }:\ 
\partial_t A_{\pm} = r A_{\pm} \mp v \partial_x A_{\pm}
 + (1 + i \alpha) \partial_x^2 A_{\pm}
 - (1 + i \beta) (\mod{A_{\pm}}^2 + \gamma \mod{A_{\mp}}^2) A_{\pm},
\nonumber 
\\
& & \phantom{1234567890}
r,v,\alpha,\beta,\gamma \hbox{ real parameters},
\nonumber 
\end{eqnarray}
Table \ref{TableResults} collects the results produced by the method
for the various PDEs in the text 
(we have added in last line an interesting system, not yet processed).
In this table,
nPq is short for ``n families of q-th order poles'',
$(m+1)^2-1$ is the maximal number of coefficients
in the subequation (\ref{eqsubeqODEOrderOnePP}),
and the column ``$a_{j,k}$'' shows the true number of such coefficients.
The column `$u_k$'' displays the minimal number of terms 
to be computed in the Laurent series (\ref{eqLaurent}).

\tabcolsep=1.5truemm
\tabcolsep=0.5truemm

\begin{table}[h] 
\caption[Summary of results.]
{Summary of results.
The solutions are labelled ``ell, trig, rat''
for, respectively:
elliptic (genus one), rational in one exponential, rational.
}
 \vspace{0.0truecm}
\begin{center}
\begin{tabular}{| l | l | l | c | c | c | l | l |}
\hline 
PDE       &poles& $m$ & $(m+1)^2-1$& $a_{j,k}$ & $u_k$ & solutions & Ref
\\ \hline 
KdV       & 1P2 &  2  &  8   &  6    &  8 & 1 ell          & 
\\ \hline 
KS        & 1P3 &  3  & 15   & 10    & 16 & 1 ell + 4 trig & \cite{MC2003}
\\ \hline 
CGL3      & 2P2 &  4  & 24   & 18    & 24 & 6 trig + 1 rat & \cite{MC2003,Hone2005}
\\ \hline 
CGL5      & 4P1 &  4  & 24   & 24    &    & 1 ell + n trig & \cite{VernovCGL5}
\\ \hline 
2 CGL3    & 2P2 &  4  & 24   & 18    &    & To be done     & \cite{CM2000b}
\\ \hline 
\end{tabular}
\end{center}
\label{TableResults}
\end{table}

\section{Possible extensions of the method}
\label{sectionExtensions}

It should first be noted that the output of the method is
all solutions which are elliptic functions of $\xi-\xi_0$
or rational in one exponential $e^{k (\xi-\xi_0)}$ 
or rational in $\xi-\xi_0$,
in which $\xi_0$ is the \textit{arbitrary} constant arising from
the first order subequation.

Having this in mind, one can think of several possible extensions.

\subsection{Extension to autonomous discrete equations}

If the independent variable $\xi$ is discrete
($\xi=n h$, with $n$ integer and $h$ a given stepsize),
can the method be extended so as to yield all elliptic
or degenerate elliptic solutions?
The question is still open.

For instance,
the discrete nonlinear equation of the Schr\"odinger type 
\begin{eqnarray}
& &
 i u_t +  p \frac{u(x+h,t)+u(x-h,t)-2 u(x)}{h^2}
 + q \frac{ {\mod u}^2 u}{1 + \nu (q h^2/p) {\mod u}^2}=0,\
p q \nu \not=0,
\nonumber \\ & & 
\label{eqdNLS_saturated}
\end{eqnarray}
admits various elliptic solitary waves \cite{KRSS,CC2008}
but no proof exists that these are the only ones.

Note that it is sufficient to look for rational functions
of $\wp(\xi),\wp'(\xi)$ since,
as shown by Briot and Bouquet \cite{BriotBouquet}, 
any elliptic function of argument $\xi+a$ can be expressed 
as a rational function of $\wp(\xi),\wp'(\xi)$.

\subsection{Extension to nonautonomous differential equations}

Consider for instance the problem 
(already solved, but this is just an example)
of finding by the above method
all the rational solutions of the second Painlev\'e equation
$\PII$
\begin{eqnarray}
& & u''=2 u^3 + x u + \alpha,
\label{eqPII}
\end{eqnarray}
i.e.~the sequence, for $\alpha$ integer,
\begin{eqnarray}
& & 
u=\mp \frac{1}{x},\ \alpha= \pm 1,
\label{eqPIIRationalSol1}
\\
& &
u=\mp \frac{2 (x^3-2)}{x(x^3+4)},\ \alpha=\pm 2,
\\
& &
u=\mp \frac{3 x^2(x^6 + 8 x^3 + 160)}{(x^3+4)(x^6+20 x^3-80)},\ \alpha=\pm 3, \dots.
\label{eqPIIRationalSol3}
\end{eqnarray} 
The integration constant $\xi_0$ of 
previous sections is no more arbitrary.

In order to retrieve these rational solutions by the method of section
\ref{sectionFirstOrderSubequation},
one must
\begin{itemize}
\item
ignore the structure of movable singularities of (\ref{eqPII});

\item
consider instead the simplest differential consequence of (\ref{eqPII})
which is autonomous, i.e.~(elimination of $x$),
\begin{eqnarray}
& & 
u u''' - u' u'' -4 u^3 u' -u^2 + \alpha u'=0;
\label{eqPIIDiffOrder3}
\end{eqnarray} 

\item
guided by (\ref{eqPIIRationalSol1})--(\ref{eqPIIRationalSol3}),
manage to find the (probably $\alpha^2$ in number) 
families of movable simple poles in (\ref{eqPIIDiffOrder3})
(it is for the moment not completely clear how this result 
should come out); 

\item
apply the method to (\ref{eqPIIDiffOrder3}) with $m=\alpha^2$
and find a unique solution.

\end{itemize}

We have not yet done it but these are at least the guidelines to be followed.

\textit{Remark}.
The symmetries of $\PII$
allow to considerably simplify the above by the change
$(u,x) \to (U,X), u=x^2 U, X=x^3$, 
thus lowering $m$ to $\mod{\alpha}$.

\subsection{Extension to Painlev\'e solutions}

After the elliptic function,
the next elementary functions defined by differential equations
are the six Painlev\'e functions, 
which all obey a nonautonomous nonlinear second order ODE.
Hence the question:
given an $N$-th order ODE (\ref{eqODE}), $N \ge 3$,
can the method of section \ref{sectionFirstOrderSubequation}
be extended so as to yield all the solutions which are algebraic transforms
of a given Painlev\'e equation?
Let us take an example to illustrate the difficulties to overcome. 

Consider the third order autonomous 
ODE built by elimination of the constant $K_1$ in
(this example is taken from the Lorenz model)
\begin{eqnarray}
& &
\frac{\D^2 x}{\D t^2}
+ 2 \frac{\D x}{\D t}
+ \frac{x^3}{2} + \left( \frac{8}{9} - \frac{K_1}{2} e^{-2 t} \right) x=0.
\label{eqLorenz2119}
\end{eqnarray}
This ODE happens to be integrable with a particular second Painlev\'e function
with $\alpha=0$ \cite{Segur}, 
\begin{eqnarray}
& &
x=a e^{-2 t/3} X,\ T=\frac{i}{2} a^{-3/2}e^{-2 t/3},\
K_1=\frac{3}{8} i a^3,
\label{eqLorenzChange}
\\
& &
\frac{\D^2 X}{\D T^2} = 2 X^3 + T X,\ 
\label{eqPIIalpha0}
\end{eqnarray}
The problem is to uncover the subequation (\ref{eqPIIalpha0}) for $X(T)$
knowing only the third order autonomous ODE for $x(t)$ and its Laurent series.
The necessity to take account of a possible change of the independent 
variable $t \to T$
introduces a countable\footnote{ 
This number can be made finite by requiring
the change $(x,t) \to (X,T)$ to introduce only exponential functions
as in (\ref{eqLorenzChange}),
so as to make the transformed ODE polynomial in the new independent
variable $T$,
before application of the method.
}
number of arbitrary coefficients 
(those of the Taylor series of $t \to T$
near the movable singularity of $x(t)$)
and thus makes it difficult to use the information provided by the
Laurent series.

In fact, the main privilege of the elliptic functions,
already mentioned as a remark in section \ref{sectionResultsBB},
is the autonomous nature of their ODE
and the invariance of this ODE under homographic transformations,
all features which do not exist any more for Painlev\'e functions.

\section{Conclusion and perspectives}
  
The main features of the method presented here are as follows:
\begin{itemize} 
\item 
It provides \textit{all} the elliptic and degenerate elliptic
(i.e.~rational in a single exponential or rational)
solutions of a given autonomous nonlinear algebraic ODE;

\item It includes all ``truncation'', ``extended'', ``new'', etc methods;

\item It makes all these methods obsolete;

\item Its main part (solving system (\ref{eqLinearSystemFj})) 
is a linear algebraic problem.

\end{itemize}

Of course, integrating the resulting first order ODE 
may consume much more time.

There exists another theory in which meromorphic functions play
a central role, this is Nevanlinna theory.
For a short comparison between Nevanlinna theory and the present method,
the interested reader can refer to \cite{ConteNg}.

Current and future work includes:

\begin{itemize}
\item Revisit physically interesting PDEs with insufficient solutions,
such as the system of two coupled CGL3 mentioned in section
\ref{sectionResults}.

\item Find general analytic solutions which are not elliptic,
for instance Lam\'e functions.
 
\end{itemize}

\section*{Acknowledgements}

Partial financial support has been provided by the
Research Grants Council contracts 
HKU 7123/05E,
HKU 7184/04E and
HKU 7038/07P.
RC acknowledges the partial support of the 
PROCORE - France/Hong Kong joint research grant
F-HK29/05T.


\vfill \eject


\begin{thebibliography}{99}

\bibitem {AblowitzClarkson}
M.~J.~Ablowitz and P.~A.~Clarkson,
\textit{Solitons, nonlinear evolution equations and inverse scattering}
(Cambridge University Press, Cambridge, 1991).

\bibitem{BriotBouquet} 
C.~Briot et J.-C.~Bouquet,
\textit{Th\'eorie des fonctions elliptiques}
(Mallet-Bachelier, Paris, 1859).
\verb+http://gallica.bnf.fr/document?O=N099571+ 

\bibitem{ChazyThese} J.~Chazy,                                
Sur les \'equations diff\'erentielles du troisi\`eme ordre et d'ordre
sup\'erieur dont l'int\'egrale g\'en\'erale a ses points critiques fixes,
Th\`ese, Paris (1910);
Acta Math.~{\bf 34} (1911) 317--385.

\bibitem{Conte1988} R.~Conte,
Universal invariance properties of Painlev\'e analysis and B\"acklund
transformation in nonlinear partial differential equations,
Phys.~Lett.~A {\bf 134} (1988) 100--104.

\bibitem{Conte1989} R.~Conte,
Invariant Painlev\'e analysis of partial differential equations,
Phys.~Lett.~A {\bf 140} (1989) 383--390.

\bibitem{CC2008} R.~Conte and K.-w.~Chow,
Doubly periodic waves of a discrete nonlinear Schr\"odinger system
with saturable nonlinearity,
J.~Nonlinear Math.~Phys.~{\bf ?} (2008) accepted.
http://arxiv.org/abs/0812.1196

\bibitem{CFP1993} R.~Conte, A.~P.~Fordy, and A.~Pickering,
A perturbative Painlev\'e approach to nonlinear differential equations,
Physica D {\bf 69} (1993) 33--58.

\bibitem{CM1989} R.~Conte and M.~Musette,
Painlev\'e analysis and B\"acklund transformation in the
Kuramoto-Sivashinsky equation,
J.~Phys.~A {\bf 22} (1989) 169--177.

\bibitem{CM1992} R.~Conte and M.~Musette,
Link between solitary waves and projective Riccati equations,
J.~Phys.~A {\bf 25} (1992) 5609--5623.

\bibitem{CM2000b} R.~Conte and M.~Musette,
Analytic expressions of hydrothermal waves,
Reports on mathematical physics {\bf 46} (2000) 77--88.
nlin.SI/0009022

\bibitem{CMBook} 
R.~Conte and M.~Musette,
\textit{The Painlev\'e handbook}, xxiv+256 pages (Springer, Berlin, 2008).
http://www.springer.com/physics/book/978-1-4020-8490-4 

\bibitem{ConteNg} R.~Conte and Ng T-w,
Meromorphic solutions of a third order nonlinear differential equation,
to be submitted (2009).

\bibitem{Darboux1883Eaux} G.~Darboux,                               
Sur les \'equations aux d\'eriv\'ees partielles,
\CRAS\ {\bf 96} (1883) 766--769.

\bibitem{FournierSpiegelThual}
J.D.~Fournier, E.A.~Spiegel and O.~Thual,   
Meromorphic integrals of two nonintegrable systems,
{\it Nonlinear dynamics},
366--373,
ed.~G.~Turchetti (World Scientific, Singapore, 1989).

\bibitem{MapleAlgcurves} 
Mark van Hoeij,
package ``algcurves'', Maple V (1997).
\hfill\break\noindent
 http://www.math.fsu.edu/\~{ }hoeij/algcurves.html

\bibitem{Hone2005} 
A.N.W.~Hone,                                  
Non-existence of elliptic travelling wave solutions of the complex
Ginzburg-Landau equation,
Physica D {\bf 205} (2005) 292--306.

\bibitem{JeffreyXu} A.~Jeffrey and Xu S.,
Travelling wave solutions to certain non-linear evolution equations,
Int.~J.~Non-Linear Mechanics {\bf 24} (1989) 425--429.

\bibitem{KRSS} A.~Khare, K.\O.~Rasmussen, M.R.~Samuelsen and A.~Saxena,
Exact solutions of the saturable discrete nonlinear Schr\"odinger equation,
J.~Phys.~A {\bf 38} (2005) 807--814.

\bibitem{KudryashovKSFourb} N.~A.~Kudryashov,
Exact soliton solutions of the generalized evolution equation
of wave dynamics,
Prikladnaia Matematika i Mekhanika {\bf 52} (1988) 465--470
[English~:
Journal of applied mathematics and mechanics {\bf 52} (1988) 361--365]

\bibitem{KudryashovElliptic} N.~A.~Kudryashov,
Exact solutions of a generalized equation of Ginzburg-Landau,
Matematicheskoye modelirovanie {\bf 1} (1989) 151--158.
[English~: none?].

\bibitem{KuramotoTsuzuki} Y.~Kuramoto and T.~Tsuzuki,
Persistent propagation of concentration waves in dissipative media far from
thermal equilibrium,
\PTP\ {\bf 55} (1976) 356--369.

\bibitem{MannevilleBook} P.~Manneville,
\textit{Dissipative structures and weak turbulence}
(Academic Press, Boston, 1990).
French adaptation:
\textit{Structures dissipatives, chaos et turbulence}
(Al\'ea-Saclay, Gif-sur-Yvette, 1991).

\bibitem{MC2003} M.~Musette and R.~Conte,
Analytic solitary waves of nonintegrable equations,
Physica D {\bf 181} (2003) 70--79.
http://arXiv.org/abs/nlin.PS/0302051

\bibitem{Pickering1993} A.~Pickering,
A new truncation in Painlev\'e analysis,
J.~Phys.~A {\bf 26} (1993) 4395--4405.

\bibitem{vS2003} W.~van Saarloos,
Front propagation into unstable states,
Physics reports {\bf 386} 29--222 (2003).

\bibitem{Segur} H.~Segur,
Solitons and the inverse scattering transform,
{\it Topics in ocean physics},
235--277,
eds.~A.~R.~Osborne and P.~Malanotte Rizzoli
(North-Holland publishing co., Amsterdam, 1982).

\bibitem{VernovCGL5} S.~Yu.~Vernov,
Elliptic solutions of the quintic complex one-dimensional Ginzburg-Landau
equation,
J.~Phys.~A {\bf 40} 9833--9844 (2007). 
\hfill\break\noindent
http://arXiv.org/abs/nlin.PS/0602060

\bibitem{WTC} J.~Weiss, M.~Tabor, and G.~Carnevale,
The Painlev\'e property for partial differential equations,
J.~Math.~Phys.~{\bf 24} (1983) 522--526.

\end{thebibliography}
\end{document}